\documentclass[a4paper]{jpconf}
\usepackage{amsmath,amsfonts,amssymb,amsthm}

\usepackage[all]{xy}

\usepackage{cite}

 \def\BAL#1\EAL{\begin{align}#1\end{align}}

 \newtheorem{theo}{Theorem}
 \newtheorem{lemma}[theo]{Lemma}
 \newtheorem{cor}[theo]{Corollary}

\DeclareMathOperator{\den}{den}

\def\R{\mathbb R}
\def\N{\mathbb N}
\def\Z{\mathbb Z}
\newcommand{\CP}{{\cal P}}
\newcommand{\CS}{{\cal S}}
 \newcommand{\nn}{\nonumber}

\begin{document}

\title{Multiplicativity in the theory of coincidence site lattices}
\author{Peter Zeiner}
\address{Fakult\"at f\"ur Mathematik, Universit\"at Bielefeld, Bielefeld, Germany}

\ead{pzeiner@math.uni-bielefeld.de}

\begin{abstract}
Coincidence Site Lattices
(CSLs) are a well established tool in the theory of grain boundaries.
For several lattices up to dimension $d=4$, the CSLs are known explicitly
as well as their indices and multiplicity functions. Many of them share
a particular property: their multiplicity functions are multiplicative.
We show how multiplicativity is connected to certain decompositions
of CSLs and the corresponding coincidence rotations and
present some criteria for multiplicativity. In general, however,
multiplicativity is violated, while supermultiplicativity still holds.
\end{abstract}

\section{Introduction}

In crystallography CSLs have been used for several decades to classify
and describe grain boundaries, see 
e.g.~\cite{boll70,gribo74,grim84,baa97} and references therein. In the
beginning research concentrated on lattices in dimensions $d\leq 3$,
but since the discovery of quasicrystals also lattices in higher dimensions
and $\Z$--modules have been studied, see e.g.~\cite{plea96,baa97,huck1}.
Here we discuss a special aspect, namely the multiplicativity of
certain combinatorial functions associated with CSLs.

Let us recall some basic concepts. 
Let
$\varGamma\subseteq\R^d$ be a $d$-dimensional lattice and
$R\in \mathrm{O}(d)$ a linear isometry. Then, $R$ is called a
(linear) \emph{coincidence isometry} of $\varGamma$
if $\varGamma(R):=\varGamma\cap R\varGamma$ is a
lattice of finite index in $\varGamma$, and
$\varGamma(R)$ is called an (ordinary or simple)
\emph{coincidence site lattice (\mbox{CSL})} (for an introduction,
see~\cite{baa97}).
The group of all coincidence isometries of $\varGamma$ is denoted by
$OC(\varGamma)$, whereas the subgroup of all orientation preserving
isometries is called the group of coincidence rotations and referred to
as $SOC(\varGamma)$.
The \emph{coincidence index} $\varSigma(R):=[\varGamma:\varGamma(R)]$
is defined as the (group theoretical) index of
$\varGamma(R)$ in $\varGamma$ and is just the ratio of the volume
of the corresponding unit cells.  

These concepts can be
generalized to include the possibility of multiple CSLs,
see~\cite{baagri05,pzmcsl2,pzmcsl3}. In particular,
the lattice
\BAL
\varGamma(R_1,\ldots,R_n):=
\varGamma\cap R_1\varGamma\cap\ldots\cap R_n\varGamma=
\varGamma(R_1)\cap\ldots\cap\varGamma(R_n)
\EAL
is called a \emph{multiple CSL (MCSL)} of order $n$, 
where $R_i$, $i\in\{1,\ldots n\}$,
are coincidence isometries of $\varGamma$. 
Its index in $\varGamma$ is denoted by
$\varSigma(R_1,\ldots,R_n)$.

There arise several interesting combinatorial questions in this
context. In particular, one is interested in the number of
coincidence isometries, coincidence rotations
and the number of CSLs of a given index $m$,
which we denote by $|\CP|f^{iso}(m)$, $|\CP'|f^{rot}(m)$ and $f(m)$,
respectively. 
Here $|\CP|$ is the order of the point group $\CP$ of $\varGamma$, and
$\CP'\subseteq \CP$ is the (normal) subgroup of orientation preserving
symmetry operations. These factors
have been chosen to guarantee that $f^{iso}(m)$ and $f^{rot}(m)$ are normalized
such that $f^{iso}(1)=1$ and $f^{rot}(1)=1$, respectively.
In general we have $f^{iso}(m)\geq f(m)$, although
$f^{iso}(m)=f^{rot}(m)=f(m)$ holds for the most important examples
in $d=2,3$ like the square lattice and the cubic lattices. In $d\geq4$
there exist several examples where $f^{iso}(m)=f^{rot}(m)>f(m)$
for infinitely many
indices $m$, see e.g. \cite{pzcsl5,pzmhliv} for examples. 

In many cases the multiplicity functions $f^{iso}(m)$,
$f^{rot}(m)$ and $f(m)$ turn out to be
multiplicative functions. Recall that a function $f:\N\to\R$
is called multiplicative if $f(mn)=f(m)f(n)$ holds
whenever $m$ and $n$ are coprime. We call $f$ supermultiplicative if
the inequality $f(mn)\geq f(m)f(n)$ holds for $m$ and $n$ coprime.
As an example we mention the square lattice~\cite{baa97}, where
\BAL
f^{iso}(m)=f^{rot}(m)=f(m)=
\begin{cases}
1& \mbox{ for } m=1\\
2^r& \mbox{ if all prime factors $p$ of $m$ satisfy $p\equiv 1\pmod 4$}\\
& \mbox{ and $r$ is the number of distinct prime factors}\\
0 & \mbox{ otherwise.}
\end{cases}
\EAL
Multiplicativity suggests to use a Dirichlet series as generating function,
which reads in this case~\cite{baa97}
\BAL
\Phi(s)
&=\sum_{m=1}^\infty\frac{f(m)}{m^s}
=\prod_{p\equiv 1(4)}\frac{1+p^{-s}}{1-p^{-s}}\nn\\
&=1+\frac{2}{5^s}+\frac{2}{13^s}+\frac{2}{17^s}+\frac{2}{25^s}+\frac{2}{29^s}
+\frac{2}{37^s}+\frac{2}{41^s}+\frac{2}{53^s}+\frac{2}{61^s}+\frac{4}{65^s}
+\frac{2}{73^s}+\ldots
\EAL
Note that it is the multiplicativity of $f(m)$ that guarantees that
the product expansion exists. In case of the square lattice, multiplicativity
is a consequence of the fact that $\Z[i]$ is a principal ideal domain.
But multiplicativity also holds for the cubic lattices in 3 dimensions
as well as for the $A_4$ root lattice and the hypercubic lattices in
4 dimensions~\cite{baa97,pzcsl4,pzcsl5}. In the latter cases multiplicativity
is due to the unique prime factorization in certain quaternion algebras.

All these cases are quite special in the sense that they are related
to algebras that allow a unique prime factorization. In fact, there
are examples where $f(m)$ and $f^{iso}(m)$ are not multiplicative, e.g.
for $\Gamma=2\Z\times 3\Z$. This raises several questions: when are
the multiplicity functions multiplicative, are there criteria for
multiplicativity? Does the multiplicativity of $f(m)$ imply
the multiplicativity of $f^{iso}(m)$ or $f^{rot}(m)$ or vice versa?
Is there a connection
between the multiplicativity of the multiplicity functions for ordinary CSLs
and multiple CSLs? What can be said if $f(m)$ is not multiplicative?

We answer some of these questions below and show that there are
some connections between multiplicativity and certain decompositions
of CSLs into multiple CSLs. We explain and motivate our results and
sketch some proofs, whereas detailed proofs will be published elsewhere.

For simplification we will consider only $f^{iso}(m)$ and $f(m)$ in the
following. In fact, this is no restriction since all properties of $f^{iso}(m)$
have a direct counterpart for $f^{rot}(m)$. Moreover $f^{rot}(m)=f^{iso}(m)$
for all lattices $\Gamma$ that contain an orientation reversing symmetry
operation in their point group, since there exists an index preserving
bijection between the symmetry preserving and symmetry reversing coincidence
isometries --- compare the remark at the beginning of Sec.~3 of \cite{pzcsl4}.

\section{Multiplicativity and supermultiplicativity of the coincidence index}

Before looking at the multiplicity functions $f(m)$, it makes sense to have
a closer look on the coincidence index $\Sigma(R)$. In particular we
are interested in $\Sigma(R_1R_2)$, if $R_1$ and $R_2$ are coincidence
isometries. We cannot expect to calculate $\Sigma(R_1R_2)$ from
$\Sigma(R_1)$ and $\Sigma(R_2)$, but we get at least
the following upper bound:
\begin{lemma}\label{lemsubSig}
$\varSigma(R_1R_2)$ divides $\varSigma(R_1)\varSigma(R_2)$.
\end{lemma}

\begin{figure}[htb]
\xymatrix@C-2em{
&\varGamma&&R_1\varGamma&&R_1R_2\varGamma\\
&{\varGamma(R_1)+\varGamma(R_1R_2)}
\ar@{-}[dl]^{n/d}\ar@{-}[dr]_k\ar@{-}[u]|*+<0.4em>{\scriptstyle m/k}
&&
{\varGamma(R_1)+R_1\varGamma(R_2)\ar@{-}[u]|*+<0.4em>{\scriptstyle d}
\ar@{-}[dl]^{m/d}\ar@{-}[dr]_{n/d}}&&
{ R_1\varGamma(R_2)+\varGamma(R_1R_2)}
\ar@{-}[dl]^k\ar@{-}[dr]_{m/d}\ar@{-}[u]|*+<0.4em>{\scriptstyle n/k}\\
\varGamma(R_1R_2)\ar@(u,l)@{-}[uur]^{mn/dk}
&&\varGamma(R_1)\ar@(u,l)@{-}[uur]^m\ar@(u,r)@{-}[uul]_m&&
R_1\varGamma(R_2)\ar@(u,r)@{-}[uul]_n\ar@(u,l)@{-}[uur]^n&&
\varGamma(R_1R_2)\ar@(u,r)@{-}[uul]_{mn/dk}\\
&&&{\hspace{-1em}\varGamma\cap R_1\varGamma\cap R_1R_2\varGamma}\hspace{-1em}
\ar@{-}[ul]_{n/d}\ar@{-}[ur]^{m/d}\ar@{-}[ulll]_k\ar@{-}[urrr]^k&
}
\caption{Relations between CSLs and their indices}
\label{figsubSig}
\end{figure}
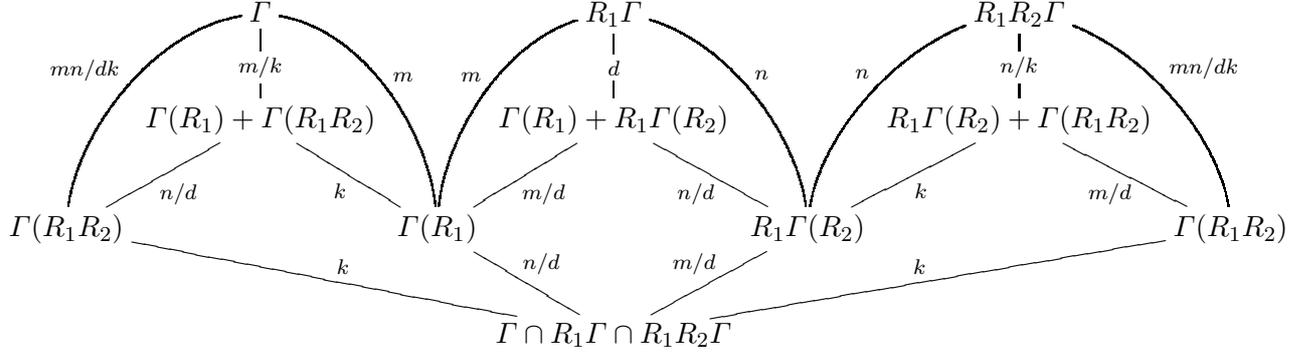

The proof makes use of the second homomorphism theorem and can be deduced
from the diagram in Fig.~\ref{figsubSig}, which shows
the relation of several CSLs and double CSLs. There we have set
$m:=\varSigma(R_1)$ and $n:=\varSigma(R_2)$. 

We can even prove more:
\begin{theo}\label{theomultSig}
If $\varSigma(R_1)$ and $\varSigma(R_2)$ are coprime, then
\BAL
\varSigma(R_1R_2)=\varSigma(R_1)\varSigma(R_2).
\EAL
\end{theo}
Note that the condition that $\varSigma(R_1)$ and $\varSigma(R_2)$ are coprime
is essential. In general we cannot expect equality. A simple counter example
is given by $R_2=R_1^{-1}$, if $\Sigma(R_1)>1$, since
$\Sigma(R_1)=\Sigma(R_1^{-1})$ holds in general~\cite{baa97}
and $\Sigma(E)=1$.

\begin{figure}[htb]
\xymatrix@C-3em{
\varGamma\hspace{1em}&&R_1\varGamma=\varGamma(R_1)+R_1\varGamma(R_2)&&
\hspace{1em}R_1R_2\varGamma\\
&\varGamma(R_1)\ar@{-}[ur]^m\ar@{-}[ul]_m&&
R_1\varGamma(R_2)\ar@{-}[ul]_n\ar@{-}[ur]^n&&\\
&\ar@(ul,dl)@{-}[uul]^{mn}&
\hspace{-2em}\varGamma(R_1R_2)
=\varGamma\cap R_1\varGamma\cap R_1R_2\varGamma\hspace{-2.4em}
\ar@{-}[ul]_n\ar@{-}[ur]^m
&\ar@(ur,dr)@{-}[uur]_{mn}
}
\caption{Relations between CSLs with $\varSigma(R_1)$ and $\varSigma(R_2)$ coprime}
\label{figmultSig}
\end{figure}
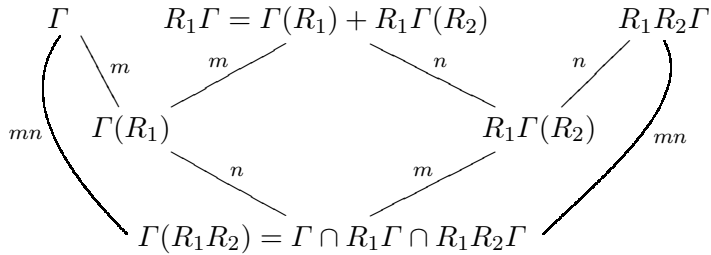

For $m$ and $n$ coprime the diagram in Fig.~\ref{figsubSig} simplifies
considerably and the result is shown in Fig.~\ref{figmultSig}. 

Moreover we can readily read off
\begin{cor}\label{cormultSig}
If $\varSigma(R_1)$ and $\varSigma(R_2)$ are coprime, then
\BAL
\varGamma(R_1R_2)=\varGamma\cap R_1\varGamma\cap R_1R_2\varGamma=
\varGamma(R_1)\cap R_1\varGamma(R_2)
\EAL
\end{cor}
This result is rather technical but plays an important role in the following,
since it relates $\varGamma(R_1R_2)$ with some kind of multiple CSLs and
provides the basis for something like a ``prime decomposition'' of CSLs.

Before we continue we want to point out that analogous results can be also
obtained for similar sublattices, see e.g. \cite{baaheu1}. 
There $\Sigma(R)$
has to be replaced by the corresponding index of the primitive
similar sublattice,
which is given by $\den(R)^d$, where $\den(R)$ is the denominator of the
similarity rotation and $d$ is the dimension. In fact, there is a close
relationship between similar sublattices and CSLs,
see~\cite{pzcsl4,svenja1}.




\section{Multiplicity functions}

We have already mentioned that the functions $f^{iso}(m)$ and $f(m)$ are
in general not multiplicative. Nevertheless  $f^{iso}(mn)$ (and likewise
$f(mn)$) is not completely independent of $f^{iso}(m)$ and $f^{iso}(n)$.
In fact we can prove~\cite{pzcslmultprep} 
\begin{theo}\label{theocrotsup}
$f^{iso}(m)$ is supermultiplicative, i.e. $f^{iso}(mn)\geq f^{iso}(m)f^{iso}(n)$
if $m$ and $n$ are coprime.
\end{theo}

The proof is mainly combinatorial and relies heavily on the multiplicativity
of the index $\Sigma$ as given in theorem~\ref{theomultSig}. 
It is this
theorem that guarantees that there are enough coincidence isometries
of index $mn$.

Recall that $|\CP|f^{iso}(m)$ counts the number of coincidence isometries
of index $m$. Correspondingly $f^{iso}(m)$ counts the number of distinct
symmetry classes $R\CP$ of coincidence isometries of index $m$, i.e.
there exist exactly $f^{iso}(m)$ different cosets $R_i \CP$ and every $R$
with $\varSigma(R)=m$ is contained in exactly one of these cosets.
Furthermore $f^{iso}(m)$ is an upper bound for the number $f(m)$ of
CSLs of index $m$. In fact, all $R\in R_i \CP$ generate the same CSL,
and hence the set of all $R$ that generate a given CSL is the union
of finitely many cosets $R_i \CP$. Let us denote the set of all 
coincidence isometries $S$ that
generate the CSL $\varGamma(R)$ by $\CS(R)$. Clearly there are exactly
$f(m)$ different $\CS(R_j)$ with index $m$. The fact that $\CS(R)$ consists
of more than one symmetry class $R\CP$ in general makes the determination
of $f(m)$ more difficult than the computation of $f^{iso}(m)$ and in general
the expressions for $f(m)$ are more complicated than those for $f^{iso}(m)$,
see~\cite{pzcsl5,pzmhliv} for examples. Nevertheless we can show
\begin{theo}\label{theocsup}
$f(m)$ is supermultiplicative, i.e. $f(mn)\geq f(m)f(n)$
if $m$ and $n$ are coprime.
\end{theo}
The proof is again combinatorial, though slightly more difficult.
In fact one needs the following lemma, which is also interesting on its own:
\begin{lemma}\label{lemcsup}
Assume that $\varSigma(R)=:m$ and $\varSigma(S)=:n$ are coprime. Then
\BAL
n\varGamma\cap\varGamma(RS)=n\varGamma(R) && \mbox{and}&&
mR\varGamma\cap\varGamma(RS)=nR\varGamma(S).
\EAL
\end{lemma}
This lemma does not only tell us that we can recover
$\varGamma(R)$ and $\varGamma(S)$ from
$\varGamma(RS)$ alone but it also tells us how to do so: just by taking the
intersection of $\varGamma(RS)$ with a suitable similar sublattice of
$\varGamma$. 

Multiplicativity can be destroyed for several reasons, and we can read them
off from this lemma.
First, there may
be isometries $Q$ of index $\varSigma(RS)=mn$ that cannot be written as
a product $Q=RS$ with $\varSigma(R)=m$ and $\varSigma(S)=n$.
As an example we mention $\Gamma=2\Z\times 3\Z$. Here $f^{iso}(6)=f(6)=1$, but
$f^{iso}(2)=f(2)=0=f^{iso}(3)=f(3)$. Further examples can be found
in~\cite{jonas}. 

Secondly, two
isometries $R$,$R'$ that generate the same CSL $\varGamma(R)=\varGamma(R')$
might give rise to different CSLs $\varGamma(RS)$ and $\varGamma(R'S)$. This
is no problem as long as $R$ and $R'$ are symmetry related. In this case
the set $\{\varGamma(R'S_k)\}_{k=1}^{f(n)}$ is just a permutation of
$\{\varGamma(RS_k)\}_{k=1}^{f(n)}$, where
$S_k$ runs over a complete set of not symmetry related $S_k$ of index
$\Sigma(S_k)=n$. However, if $R$ and $R'$ are not
symmetry related additional CSLs might occur. 

Analogous results hold for similar sublattices. Let $g(m)$ be the number of
similar sublattices of index $m$. 
Then the function $g(m)$ is in general only
supermultiplicative. An example
for a lattice with non--multiplicative $g(m)$ is again
the lattice $\Gamma=2\Z\times 3\Z$. But note that similar sublattices seem
to be more sensitive to violation of multiplicativity than CSLs.
E.g., for $\Gamma=\Z\times 5\Z$ multiplicativity is violated for $g(m)$
while $f(m)$ is still multiplicative~\cite{pzssl1,jonas}.

\section{A criterion for multiplicativity}

We have seen that $f(m)$ and $f^{iso}(m)$ are in general only
supermultiplicative. Now the interesting question is whether there exist
some criteria for multiplicativity and the answer is positive.
A first hint is given by known examples
in $d\leq 4$. For root lattices in $d\leq 4$ the multiplicity functions
$f(m)$ and $f^{iso}(m)$
are usually multiplicative. The reason is that these lattices are related 
to principal ideal domains (and thus unique factorization domains)
of algebraic integers or quaternions. So we expect
that some kind of unique factorization property is essential.
In fact we can prove
\begin{theo}\label{mc}
The following statements are equivalent:
\begin{enumerate}
\item \label{mc1} $f(m)$ is multiplicative.
\item \label{mc2} 
Every (ordinary) CSL $\varGamma(R)$ can be written (uniquely) as
$\varGamma(R)=\varGamma(R_1)\cap\ldots\cap\varGamma(R_n)$, where the indices
$\varSigma(R_i)$ are powers of distinct primes.
\item \label{mc3} 
Every MCSL $\varGamma(R_1,\ldots,R_n)$ of order $n$
can be written (uniquely) as
$\varGamma(R_1,\ldots,R_n)=\varGamma_1\cap\ldots\cap\varGamma_k$,
where the $\varGamma_k$
are MCSLs of order at most $n$ and whose indices
$\varSigma_k$ are powers of distinct primes.
\end{enumerate}
\end{theo}
Note that Lemma~\ref{lemcsup} guarantees the uniqueness of
the decomposition $\varGamma(R)=\varGamma(R_1)\cap\ldots\cap\varGamma(R_n)$,
if it exists. 
Lemma~\ref{lemcsup} is also the main ingredient for the proof, which again
also involves combinatorial arguments. 
But be careful. The decomposition
$\varGamma(R)=\varGamma(R_1)\cap\ldots\cap\varGamma(R_n)$ does not imply
the decomposition $R=R_1\cdots R_n$, in general $R$ is not even symmetry
related to $R_1\cdots R_n$. The situation is even worse. If
$R=R_1\cdots R_n$ is a decomposition of $R$ then in general
$\varGamma(R)\neq\varGamma(R_1)\cap\ldots\cap\varGamma(R_n)$, which is due
to the fact that $O(n)$ is not Abelian for $n\geq 2$.

A corresponding criterion for the coincidence isometries exists as well.
The formulation of it is a bit more intricate, since isometries usually do not
commute. For CSLs the decomposition into its prime power constituents is
unique (up to permutation), for isometries a decomposition will depend strongly
on how the factors are ordered. 

First notice that if the coincidence isometry $R$ with $\varSigma(R)=mn$
can be factored as $R=R_1R_2$ with $\varSigma(R_1)=m$ and
$\varSigma(R_2)=n$ coprime, then $R_1$ and $R_2$ are uniquely determined up
to elements of the point group $\CP$, i.e. all other decompositions
are of the form  $R=(R_1Q)(Q^{-1}R_2)$ with $Q\in\CP$. 
Note that
$R_2$ and $Q^{-1}R_2$ are usually not symmetry related, whereas $R_1$ and
$R_1Q$ are. 

At this point it is not clear whether the existence of
a decomposition $R=R_1R_2$ implies a decomposition $R=R'_2R'_1$, where
$\varSigma(R_1)=\varSigma(R'_1)=m$ and $\varSigma(R_2)=\varSigma(R'_2)=n$.
This motivates the following definitions:
We call a bijection $\pi=\{p_1,p_2\ldots\}$
from the positive integers onto the prime numbers
an ordering of the prime numbers. We call a decomposition of a coincidence
isometry $R=R_1\cdots R_n$ a $\pi$--decomposition of $R$ if
$\varSigma(R_i)$ is a power of $p_i$ for any $i$
(we allow $\varSigma(R_i)=p_i^0=1$).
It is clear that any 
$\pi$--decomposition can be unique only up to point group elements.

We can now formulate the analogue of theorem~\ref{mc} for $f^{iso}(m)$
\begin{theo}
The following statements are equivalent:
\begin{enumerate}
\item \label{mcrot1} $f^{iso}(m)$ is multiplicative.
\item \label{mcrot2} There exists an ordering $\pi$ of the prime numbers
such that any coincidence isometry $R$ has a (unique) $\pi$--decomposition.
\item \label{mcrot3} For any ordering $\pi$ of the prime numbers
there exists a $\pi$--decomposition of every coincidence isometry $R$.
\end{enumerate}
\end{theo}

Given these two quite similar criteria we may expect that there
is some connection between the multiplicativity of $f(m)$ and $f^{iso}(m)$.
In fact we can prove
\begin{theo}\label{iso2csl}
$f(m)$ is multiplicative if $f^{iso}(m)$ is.
\end{theo}

Here it is worth to comment on the various decompositions that occur here.
For simplicity we assume that only two prime powers are involved,
say $\Sigma(R)=p_1^{r_1}p_2^{r_2}$. Then the multiplicativity of $f^{iso}(m)$
guarantees the existence of two decompositions $R=R_1S_1$ and $R=R_2S_2$
with $\Sigma(R_1)=p_1^{r_1}=\Sigma(S_2)$ and $\Sigma(R_2)=p_2^{r_2}=\Sigma(S_1)$.
In this case the unique decomposition of $\Gamma(R)$ reads
$\varGamma(R)=\varGamma(R_1)\cap\varGamma(R_2)$. So given the decompositions of
$R$ we get immediately the decomposition of $\varGamma(R)$. However, it does
not work the other way round. So given a decomposition of $\varGamma(R)$
we do not get any information on the decompositions of $R$, even if we
would know that they exist.

Now what about the converse of theorem~\ref{iso2csl}? Does it exist? 
We do not know the answer so far. So we conclude with an open question:
Does the multiplicativity of $f(m)$ imply the multiplicativity of
$f^{iso}(m)$? Or are there lattices with
multiplicative $f(m)$ but non--multiplicative $f^{iso}(m)$?

\section*{Acknowledgements}

This work 
was supported by the German Research Council (DFG), within the CRC~701.

\section*{References}
\bibliography{liv1}

\end{document}